\documentclass[a4paper,11pt]{article}

\usepackage{amsmath}
\usepackage{amsfonts}
\usepackage{amssymb}
\usepackage{amsthm}

\newtheorem{theorem}{Theorem}[section]
\newtheorem{lemma}[theorem]{Lemma}
\newtheorem{proposition}[theorem]{Proposition}

\theoremstyle{definition}
\newtheorem{definition}[theorem]{Definition}

\theoremstyle{remark}

\numberwithin{equation}{section}

\begin{document}

\title{Uniqueness for a Stochastic Inviscid\\ Dyadic Model}

\author{D. Barbato, F. Flandoli and F. Morandin}


\date{October 21st, 2009}


\maketitle

\begin{abstract}
For the deterministic dyadic model of turbulence, there are examples
of initial conditions in $l^{2}$ which have more than one
solution. The aim of this paper is to prove that uniqueness, for all
$l^{2}$-initial conditions, is restored when a suitable multiplicative
noise is introduced. The noise is formally energy
preserving. Uniqueness is understood in the weak probabilistic sense.
\end{abstract}


\section{Introduction}

The infinite system of nonlinear differential equations
\begin{align}
\frac{dX_{n}\left(  t\right)  }{dt}  &  =k_{n-1}X_{n-1}^{2}\left(  t\right)
-k_{n}X_{n}\left(  t\right)  X_{n+1}\left(  t\right)  ,\quad t\geq
0\label{det dyadic}\\
X_{n}\left(  0\right)   &  =x_{n}\nonumber
\end{align}
for $n\geq1$, with coefficients $k_{n}>0$ for each $n\geq1$, $X_{0}(t)=0$ and
$k_{0}=0$, is one of the simplest models which presumably reflect some of the
properties of 3D Euler equations. At least, it is infinite dimensional,
formally conservative (the energy $\sum_{n=1}^{\infty}X_{n}^{2}\left(
t\right)  $ is formally constant), and quadratic. One of its `pathologies' is
the lack of uniqueness of solutions, in the space $l^{2}$ of square summable
sequences: when, for instance, $k_{n}=\lambda^{n}$ with $\lambda>1$, there are
examples of initial conditions $x=\left(  x_{n}\right)  _{n\geq1}\in l^{2}$
such that there exists at least two solutions in $l^{2}$ on some interval
$\left[  0,T\right]  $, with continuous components. This has been proved in
\cite{BFM}: on one side, given any $x\in l^{2}$, there exists a solution such
that $\sum_{n=1}^{\infty}X_{n}^{2}\left(  t\right)  \leq\sum_{n=1}^{\infty
}x_{n}^{2}$ for all $t\geq0$; on the other side, for special elements
$a=\left(  a_{n}\right)  _{n\geq1}\in l^{2}$, the strictly increasing
sequence
\[
\frac{a_{n}}{t_{0}-t},\quad t\in\lbrack0,t_{0})
\]
is a (selfsimilar) solution. Other counterexamples can be done by
time-reversing any solution which dissipate energy (this happens for all
solutions having positive components). System (\ref{det dyadic}) and variants
of it have other special features, like energy dissipation and loss of
regularity, see \cite{Ces}, \cite{CFP1}, \cite{CFP2}, \cite{FriPav},
\cite{KatPav}, \cite{KiZlat}, \cite{Wal}, \cite{BFM}.

In this paper we prove that uniqueness is restored under a suitable random
perturbation. On a filtered probability space $\left(  \Omega,F_{t},P\right)
$, let $\left(  W_{n}\right)  _{n\geq1}$ be a sequence of independent Brownian
motions. We consider the infinite system of stochastic differential equations
in Stratonovich form
\begin{equation}
dX_{n}=\left(  k_{n-1}X_{n-1}^{2}-k_{n}X_{n}X_{n+1}\right)  dt+\sigma
k_{n-1}X_{n-1}\circ dW_{n-1}-\sigma k_{n}X_{n+1}\circ dW_{n}
\label{eq stratonovich}%
\end{equation}
for $n\geq1$, with $X_{0}(t)=0$ and $\sigma\neq0$. The concept of
exponentially integrable solution, used in the following theorem, is defined
in the next section. By classical arguments, we shall prove weak existence in
the class of exponentially integrable solution. Our main theorem is:

\begin{theorem}
\label{main uniqueness thm}Given $x\in l^{2}$, in the class of exponentially
integrable solutions on an interval $\left[  0,T\right]  $ there is weak
uniqueness for equation (\ref{eq stratonovich}).
\end{theorem}

The proof is given in section \ref{sect Girsanov}. Weak uniqueness here means
uniqueness of the law of the process on the space $C\left(  \left[
0,T\right]  ;\mathbb{R}\right)  ^{\mathbb{N}}$. Our approach is based on
Girsanov transformation, so this is the natural result one expects. We do not
know about strong uniqueness. Our use of Girsanov transformation is not the
most classical one, and is inspired to \cite{AirMall}, \cite{CruzFM}.

The multiplicative noise in equation (\ref{eq stratonovich}) preserves the
formal energy conservation. By applying the rules of Stratonovich calculus
(see the same computation at the It\^{o} level in the proof of Theorem
\ref{teor uniq linear}) we have%
\begin{align*}
dX_{n}^{2}  &  =2\left(  k_{n-1}X_{n-1}^{2}X_{n}-k_{n}X_{n}^{2}X_{n+1}\right)
dt\\
&  +\sigma k_{n-1}X_{n-1}X_{n}\circ dW_{n-1}-\sigma k_{n}X_{n}X_{n+1}\circ
dW_{n}%
\end{align*}
so that, \textit{formally}, $d\sum_{n=1}^{\infty}X_{n}^{2}\left(  t\right)
=0$. Only a multiplicative noise of special form has this property, which is
one of the key properties formally verified by Euler equations. Notice that
the It\^{o} formulation (\ref{eq Ito}) below (see also the linear analog
(\ref{linear Ito})), contains a dissipative term, which however is exactly
balanced by the correction term when It\^{o} formula is applied. Thus equation
(\ref{eq Ito}) below is not formally dissipative as it may appear at first glance.

Having in mind the lack of uniqueness, or at least the open problems about
uniqueness, typical of various deterministic models in fluid dynamics, we
think it is relevant to know that suitable stochastic perturbations may
restore uniqueness. An example in this direction is know for the linear
transport equation with poor regularity of coefficients, see \cite{FGP}. The
model of the present paper seems to be the first nonlinear example of this
regularization phenomenon (in the area of equations of fluid dynamic type,
otherwise see \cite{Gyongy}, \cite{GyongyPard} and related works, based on
completely different methods). Partial results in the direction of
improvements of well posedness, by means of additive noise, have been obtained
for the 3D Navier-Stokes equations and other models by \cite{DaPDeb}, \cite{FR
Markov}.

Let us remark that, although equation (\ref{eq stratonovich}) is not a PDE, it
has a vague correspondence with the stochastic Euler equation
\[
du+\left[  u\cdot\nabla u+\nabla p\right]  dt+\sigma\sum_{j}\nabla u\circ
dW^{j}\left(  t\right)  =0,\quad\mathrm{div}u=0.
\]
The energy is formally conserved also in this equation. The noise of this
equation is multiplicative as in \cite{FGP}, linearly dependent on first
derivatives of the solution. For a Lagrangian motivation of such a noise, in
the case of stochastic Navier-Stokes equations, see \cite{Mik Roz}.

In equation (\ref{eq stratonovich}), we have inserted the parameter
$\sigma\neq0$ just to emphasize the basic open problem of understanding the
zero-noise limit, $\sigma\rightarrow0$. For simple examples of linear
transport equations this is possible and yields a nontrivial selection
principle among different solutions of the deterministic limit equation, see
\cite{AttFla}. In the nonlinear case of the present paper the small
coefficient $\sigma$ appears in the form of a singular perturbation in the
Girsanov density, thus the analysis of $\sigma\rightarrow0$ is nontrivial.

\section{It\^{o} formulation}

The It\^{o} form of equation (\ref{eq stratonovich}) is
\begin{align}
dX_{n}  &  =\left(  k_{n-1}X_{n-1}^{2}-k_{n}X_{n}X_{n+1}\right)  dt+\sigma
k_{n-1}X_{n-1}dW_{n-1}-\sigma k_{n}X_{n+1}dW_{n}\label{eq Ito}\\
&  -\frac{\sigma^{2}}{2}\left(  k_{n}^{2}+k_{n-1}^{2}\right)  X_{n}dt\nonumber
\end{align}
for all $n\geq1$, with $k_{0}=0$ and $X_{0}=0$, as explained at the end of
this section. All our rigorous analyses are based on the It\^{o} form, the
Stratonovich one serving mainly as an heuristic guideline.

Let us introduce the concept of weak solution (equivalent to the concept of
solution of the martingale problem). Since our main emphasis is on uniqueness,
we shall always restrict ourselves to a finite time horizon $\left[
0,T\right]  $.

By a filtered probability space $\left(  \Omega,F_{t},P\right)  $, on a finite
time horizon $\left[  0,T\right]  $, we mean a probability space $\left(
\Omega,F_{T},P\right)  $ and a right-continuous filtration $\left(
F_{t}\right)  _{t\in\left[  0,T\right]  }$.

\begin{definition}
\label{def sol 1}Given $x\in l^{2}$, a weak solution of equation (\ref{eq
stratonovich}) in $l^{2}$ is a filtered probability space $\left(
\Omega,F_{t},P\right)  $, a sequence of independent Brownian motions $\left(
W_{n}\right)  _{n\geq1}$ on $\left(  \Omega,F_{t},P\right)  $, and an $l^{2}%
$-valued stochastic process $\left(  X_{n}\right)  _{n\geq1}$ on $\left(
\Omega,F_{t},P\right)  $, with continuous adapted components $X_{n}$, such
that
\begin{align*}
X_{n}\left(  t\right)   &  =x_{n}+\int_{0}^{t}\left(  k_{n-1}X_{n-1}%
^{2}\left(  s\right)  -k_{n}X_{n}\left(  s\right)  X_{n+1}\left(  s\right)
\right)  ds\\
&  +\int_{0}^{t}\sigma k_{n-1}X_{n-1}\left(  s\right)  dW_{n-1}\left(
s\right)  -\int_{0}^{t}\sigma k_{n}X_{n+1}\left(  s\right)  dW_{n}\left(
s\right)  \\
&  -\int_{0}^{t}\frac{\sigma^{2}}{2}\left(  k_{n}^{2}+k_{n-1}^{2}\right)
X_{n}\left(  s\right)  ds
\end{align*}
for each $n\geq1$, with $k_{0}=0$ and $X_{0}=0$. We denote this solution by
$\left(  \Omega,F_{t},P,W,X\right)  $, or simply by $X$.
\end{definition}

To prove uniqueness we need the following technical condition, that we call
exponential integrability, for shortness.

\begin{definition}
We say that a weak solution $\left(  \Omega,F_{t},P,W,X\right)  $ is
exponentially integrable if
\[
E^{P}\left[  e^{\frac{1}{\sigma^{2}}\int_{0}^{T}\sum_{n=1}^{\infty}X_{n}%
^{2}\left(  t\right)  dt}\left(  1+\int_{0}^{T}X_{i}^{4}\left(  t\right)
dt\right)  ^{2}\right]  <\infty
\]
for all $i\in\mathbb{N}$.

We say that a weak solution is of class $L^{\infty}$ if there is a constant
$C>0$ such that $\sum_{n=1}^{\infty}X_{n}^{2}\left(  t\right)  \leq C$ for
a.e. $\left(  \omega,t\right)  \in\Omega\times\left[  0,T\right]  $.
\end{definition}

$L^{\infty}$-solutions are exponentially integrable. Our main result, theorem
\ref{main uniqueness thm}, states the weak uniqueness in the class of
exponentially integrable solutions. In addition, we have:

\begin{theorem}
\label{teo existence}Given $\left(  x_{n}\right)  \in l^{2}$, there exits a
weak $L^{\infty}$-solution to equation (\ref{eq stratonovich}).
\end{theorem}

The proof is given in section \ref{sect Girsanov} and is based again on
Girsanov transform. However, we remark that existence can be proved also by
compactness method, similarly to the case of stochastic Euler or Navier-Stokes
equations, see for instance \cite{BrePes} and \cite{FlCime}. In both cases,
notice that it is a weak existence result: the solution is not necessarily
adapted to the completed filtration of the Brownian motions.

The following proposition clarifies that a process satisfying (\ref{eq Ito})
rigorously satisfies also (\ref{eq stratonovich}).

\begin{proposition}
If $X$ is a weak solution of equation (\ref{eq stratonovich}), then for every
$n\geq1$ the process $\left(  X_{n}\left(  t\right)  \right)  _{t\geq0}$ is a
continuous semimartingale, hence the two Stratonovich integrals
\begin{align*}
&  \int_{0}^{t}k_{n-1}X_{n-1}\left(  s\right)  \circ dW_{n-1}\left(  s\right)
\quad\text{ for }n\geq2\\
&  -\int_{0}^{t}k_{n}X_{n+1}\left(  s\right)  \circ dW_{n}\left(  s\right)
\quad\text{ for }n\geq1
\end{align*}
are well defined and equal, respectively, to
\begin{align*}
&  \int_{0}^{t}k_{n-1}X_{n-1}\left(  s\right)  dW_{n-1}\left(  s\right)
-\frac{\sigma}{2}\int_{0}^{t}k_{n-1}^{2}X_{n}\left(  s\right)  ds\\
&  -\int_{0}^{t}k_{n}X_{n+1}\left(  s\right)  dW_{n}\left(  s\right)
-\frac{\sigma}{2}\int_{0}^{t}k_{n}^{2}X_{n}\left(  s\right)  ds.
\end{align*}
Hence $X$ satisfies the Stratonovich equations (\ref{eq stratonovich}).
\end{proposition}

\begin{proof}
We use a number of concepts and rules of stochastic calculus that can be found
for instance in \cite{Kun}. We have
\[
\int_{0}^{t}X_{n-1}\left(  s\right)  \circ dW_{n-1}\left(  s\right)  =\int
_{0}^{t}X_{n-1}\left(  s\right)  dW_{n-1}\left(  s\right)  +\frac{1}{2}\left[
X_{n-1},W_{n-1}\right]  _{t}%
\]
where $\left[  X_{n-1},W_{n-1}\right]  _{t}$ is the joint quadratic variation
of $X_{n-1}$ and $W_{n-1}$. From the equation for $X_{n-1}\left(  t\right)  $,
using the independence of the Brownian motions, we can compute $\left[
X_{n-1},W_{n-1}\right]  _{t}=-\int_{0}^{t}\sigma k_{n-1}X_{n}\left(  s\right)
ds$. Similarly
\[
\int_{0}^{t}X_{n+1}\left(  s\right)  \circ dW_{n}\left(  s\right)  =\int
_{0}^{t}X_{n+1}\left(  s\right)  dW_{n}\left(  s\right)  +\frac{1}{2}\left[
X_{n+1},W_{n}\right]  _{t}%
\]
and $\left[  X_{n+1},W_{n}\right]  _{t}=\int_{0}^{t}\sigma k_{n}X_{n}\left(
s\right)  ds$. The proof is complete.
\end{proof}

\section{Auxiliary linear equation}

Up to Girsanov transform (section \ref{sect Girsanov}), our results are based
on the following infinite system of \textit{linear} stochastic differential
equations
\begin{align*}
dX_{n}  &  =\sigma k_{n-1}X_{n-1}\circ dB_{n-1}-\sigma k_{n}X_{n+1}\circ
dB_{n}\\
X_{n}\left(  0\right)   &  =x_{n}%
\end{align*}
for $n\geq1$, with $X_{0}(t)=0$ and $\sigma\neq0$, where $\left(
B_{n}\right)  _{n\geq0}$ is a sequence of independent Brownian motions. The
It\^{o} formulation is
\begin{align}
dX_{n}  &  =\sigma k_{n-1}X_{n-1}dB_{n-1}-\sigma k_{n}X_{n+1}dB_{n}%
-\frac{\sigma^{2}}{2}\left(  k_{n}^{2}+k_{n-1}^{2}\right)  X_{n}%
dt\label{linear Ito}\\
X_{n}\left(  0\right)   &  =x_{n}.\nonumber
\end{align}

\begin{definition}
Let $\left(  \Omega,F_{t},Q\right)  $ be a filtered probability space and let
$\left(  B_{n}\right)  _{n\geq0}$ be a sequence of independent Brownian
motions on $\left(  \Omega,F_{t},Q\right)  $. Given $x\in l^{2}$, a solution
of equation (\ref{linear Ito}) on $\left[  0,T\right]  $ in the space $l^{2}$
is an $l^{2}$-valued stochastic process $\left(  X\left(  t\right)  \right)
_{t\in\left[  0,T\right]  }$, with continuous adapted components $X_{n}$, such
that $Q$-a.s.
\begin{align*}
X_{n}\left(  t\right)   &  =x_{n}+\int_{0}^{t}\sigma k_{n-1}X_{n-1}\left(
s\right)  dB_{n-1}\left(  s\right)  -\int_{0}^{t}\sigma k_{n}X_{n+1}\left(
s\right)  dB_{n}\left(  s\right)  \\
&  -\int_{0}^{t}\frac{\sigma^{2}}{2}\left(  k_{n}^{2}+k_{n-1}^{2}\right)
X_{n}\left(  s\right)  ds
\end{align*}
for each $n\geq1$ and $t\in\left[  0,T\right]  $, with $k_{0}=0$ and $X_{0}=0$.
\end{definition}

Our main technical result is the following theorem.

\begin{theorem}
\label{teor uniq linear}Given $x\in l^{2}$, in the class of solutions of
equation (\ref{linear Ito}) on $\left[  0,T\right]  $ such that
\begin{equation}
\int_{0}^{T}E^{Q}\left[  X_{n}^{4}\left(  t\right)  \right]  dt<\infty
\label{forth order}%
\end{equation}
for each $n\geq1$ and
\begin{equation}
\lim_{n\rightarrow\infty}\int_{0}^{T}E^{Q}\left[  X_{n}^{2}\left(  t\right)
\right]  dt=0\label{assumption limit}%
\end{equation}
there is at most one element.
\end{theorem}

\begin{proof}
By linearity, it is sufficient to prove that a solution $\left(  X_{n}\right)
_{n\geq1}$, with properties (\ref{forth order}) and (\ref{assumption limit}),
with null initial condition is the zero solution. Assume thus $x=0$. We have
\begin{align*}
X_{n}\left(  t\right)   &  =x_{n}+\int_{0}^{t}\sigma k_{n-1}X_{n-1}\left(
s\right)  dB_{n-1}\left(  s\right)  -\int_{0}^{t}\sigma k_{n}X_{n+1}\left(
s\right)  dB_{n}\left(  s\right)  \\
&  -\int_{0}^{t}\frac{\sigma^{2}}{2}\left(  k_{n}^{2}+k_{n-1}^{2}\right)
X_{n}\left(  s\right)  ds
\end{align*}
hence, from It\^{o} formula, we have
\begin{align*}
\frac{1}{2}dX_{n}^{2} &  =X_{n}dX_{n}+\frac{1}{2}d\left[  X_{n}\right]  _{t}\\
&  =-\frac{\sigma^{2}}{2}\left(  k_{n}^{2}+k_{n-1}^{2}\right)  X_{n}%
^{2}dt+dM_{n}+\frac{\sigma^{2}}{2}\left(  k_{n-1}^{2}X_{n-1}^{2}+k_{n}%
^{2}X_{n+1}^{2}\right)  dt
\end{align*}
where
\[
M_{n}\left(  t\right)  =\int_{0}^{t}\sigma k_{n-1}X_{n-1}\left(  s\right)
X_{n}\left(  s\right)  dB_{n-1}\left(  s\right)  -\int_{0}^{t}\sigma
k_{n}X_{n}\left(  s\right)  X_{n+1}\left(  s\right)  dB_{n}\left(  s\right)  .
\]
From (\ref{forth order}), $M_{n}\left(  t\right)  $ is a martingale, for each
$n\geq1$, hence $E^{Q}\left[  M_{n}\left(  t\right)  \right]  =0$. Moreover,
for each $n\geq1$, $E^{Q}\left[  X_{n}^{2}\left(  t\right)  \right]  $ is
finite and continuous in $t$: it follows easily from condition (\ref{forth
order}) and equation (\ref{linear Ito}) itself. From the previous equation
(and the property $E^{Q}\left[  X_{n}^{2}\left(  0\right)  \right]  =0)$ we
deduce that $E^{Q}\left[  X_{n}^{2}\left(  t\right)  \right]  $ satisfies
\begin{align*}
E^{Q}\left[  X_{n}^{2}\left(  t\right)  \right]   &  =-\sigma^{2}\left(
k_{n}^{2}+k_{n-1}^{2}\right)  \int_{0}^{t}E^{Q}\left[  X_{n}^{2}\left(
s\right)  \right]  ds\\
&  +\sigma^{2}k_{n-1}^{2}\int_{0}^{t}E^{Q}\left[  X_{n-1}^{2}\left(  s\right)
\right]  ds+\sigma^{2}k_{n}^{2}\int_{0}^{t}E^{Q}\left[  X_{n+1}^{2}\left(
s\right)  \right]  ds
\end{align*}
for $n\geq1$, with $u_{0}\left(  t\right)  =0$ for $t\geq0$. It follows
\[
\int_{0}^{t}E^{Q}\left[  \left(  X_{n+1}^{2}\left(  s\right)  -X_{n}%
^{2}\left(  s\right)  \right)  \right]  ds\geq\frac{k_{n-1}^{2}}{k_{n}^{2}%
}\int_{0}^{t}E^{Q}\left[  \left(  X_{n}^{2}\left(  s\right)  -X_{n-1}%
^{2}\left(  s\right)  \right)  \right]  ds.
\]
Since $X_{0}\equiv0$, we have $\int_{0}^{t}E^{Q}\left[  \left(  X_{1}%
^{2}\left(  s\right)  -X_{0}^{2}\left(  s\right)  \right)  \right]  ds\geq0$
and thus
\[
\int_{0}^{t}E^{Q}\left[  \left(  X_{n+1}^{2}\left(  s\right)  -X_{n}%
^{2}\left(  s\right)  \right)  \right]  ds\geq0
\]
for every $n\geq1$, by induction. This implies
\[
\int_{0}^{T}E^{Q}\left[  X_{n}^{2}\left(  s\right)  \right]  ds\leq\int
_{0}^{T}E^{Q}\left[  X_{n+1}^{2}\left(  s\right)  \right]  ds
\]
for all $n\geq1$. Therefore, by assumption (\ref{assumption limit}), for every
$n\geq1$ we have $\int_{0}^{T}E^{Q}\left[  X_{n}^{2}\left(  s\right)  \right]
ds=0$. This implies $X_{n}^{2}\left(  s\right)  =0$ a.s. in $\left(
\omega,s\right)  $, hence $X$ is the null process. The proof is complete.
\end{proof}

We complete this section with an existence result. The class $L^{\infty
}\left(  \Omega\times\left[  0,T\right]  ;l^{2}\right)  $ is included in the
class described by the uniqueness theorem.

Notice that this is a result of strong existence and strong (or pathwise) uniqueness.

\begin{theorem}
\label{teo existence linear}Given $x\in l^{2}$ , there exists a unique
solution in $L^{\infty}\left(  \Omega\times\left[  0,T\right]  ;l^{2}\right)
$, with continuous components.
\end{theorem}

\begin{proof}
We have only to prove existence. For every positive integer $N$, consider the
finite dimensional stochastic system
\begin{align*}
dX_{n}^{\left(  N\right)  }  &  =\sigma k_{n-1}X_{n-1}^{\left(  N\right)
}dB_{n-1}-\sigma k_{n}X_{n+1}^{\left(  N\right)  }dB_{n}-\frac{\sigma^{2}}%
{2}\left(  k_{n}^{2}+k_{n-1}^{2}\right)  X_{n}^{\left(  N\right)  }dt\\
X_{n}^{\left(  N\right)  }\left(  0\right)   &  =x_{n}%
\end{align*}
for $n=1,...,N$, with $k_{0}=k_{N}=0$, $X_{0}^{\left(  N\right)  }\left(
t\right)  =X_{N+1}^{\left(  N\right)  }\left(  t\right)  =0$. This linear
finite dimensional equation has a unique global strong solution. By It\^{o}
formula
\[
\frac{1}{2}d\left(  X_{n}^{\left(  N\right)  }\right)  ^{2}=X_{n}^{\left(
N\right)  }dX_{n}^{\left(  N\right)  }+\frac{1}{2}d\left[  X_{n}^{\left(
N\right)  },X_{n}^{\left(  N\right)  }\right]  _{t}%
\]%
\begin{align*}
&  =\sigma k_{n-1}X_{n}^{\left(  N\right)  }X_{n-1}^{\left(  N\right)
}dB_{n-1}-\sigma k_{n}X_{n}^{\left(  N\right)  }X_{n+1}^{\left(  N\right)
}dB_{n}-\frac{\sigma^{2}}{2}\left(  k_{n}^{2}+k_{n-1}^{2}\right)  \left(
X_{n}^{\left(  N\right)  }\right)  ^{2}dt\\
&  +\frac{\sigma^{2}}{2}\left(  k_{n-1}^{2}\left(  X_{n-1}^{\left(  N\right)
}\right)  ^{2}+k_{n}^{2}\left(  X_{n+1}^{\left(  N\right)  }\right)
^{2}\right)
\end{align*}
hence
\begin{align*}
\frac{1}{2}d\sum_{n=1}^{N}\left(  X_{n}^{\left(  N\right)  }\right)  ^{2}  &
=\sum_{n=1}^{N}\sigma k_{n-1}X_{n}^{\left(  N\right)  }X_{n-1}^{\left(
N\right)  }dB_{n-1}-\sum_{n=1}^{N}\sigma k_{n}X_{n}^{\left(  N\right)
}X_{n+1}^{\left(  N\right)  }dB_{n}\\
&  -\frac{\sigma^{2}}{2}\sum_{n=1}^{N}k_{n}^{2}\left(  X_{n}^{\left(
N\right)  }\right)  ^{2}dt+\frac{\sigma^{2}}{2}\sum_{n=1}^{N}k_{n-1}%
^{2}\left(  X_{n-1}^{\left(  N\right)  }\right)  ^{2}\\
&  -\frac{\sigma^{2}}{2}\sum_{n=1}^{N}k_{n-1}^{2}\left(  X_{n}^{\left(
N\right)  }\right)  ^{2}dt+\frac{\sigma^{2}}{2}\sum_{n=1}^{N}k_{n}^{2}\left(
X_{n+1}^{\left(  N\right)  }\right)  ^{2}.
\end{align*}
This is equal to zero. Thus
\[
\sum_{n=1}^{N}\left(  X_{n}^{\left(  N\right)  }\right)  ^{2}\left(  t\right)
=\sum_{n=1}^{N}x_{n}^{2},\quad Q\text{-a.s.}%
\]
In particular, this very strong bound implies that there exists a subsequence
$N_{k}\rightarrow\infty$ such that $\left(  X_{n}^{\left(  N_{k}\right)
}\right)  _{n\geq1}$ converges weakly to some $\left(  X_{n}\right)  _{n\geq
1}$ in $L^{p}\left(  \Omega\times\left[  0,T\right]  ;l^{2}\right)  $ for
every $p>1$ and also weak star in $L^{\infty}\left(  \Omega\times\left[
0,T\right]  ;l^{2}\right)  $. Hence in particular $\left(  X_{n}\right)
_{n\geq1}$ belongs to $L^{\infty}\left(  \Omega\times\left[  0,T\right]
;l^{2}\right)  $. Now the proof proceeds by standard arguments typical of
equations with monotone operators (which thus apply to linear equations),
presented in \cite{Pardoux}, \cite{KrylovRoz}, The subspace of $L^{p}\left(
\Omega\times\left[  0,T\right]  ;l^{2}\right)  $ of progressively measurable
processes is strongly closed, hence weakly closed, hence $\left(
X_{n}\right)  _{n\geq1}$ is progressively measurable. The one-dimensional
stochastic integrals which appear in each equation of system (\ref{linear
Ito}) are (strongly) continuous linear operators from the subspace of
$L^{2}\left(  \Omega\times\left[  0,T\right]  ;l^{2}\right)  $ of
progressively measurable processes to $L^{2}\left(  \Omega\right)  $, hence
they are weakly continuous, a fact that allows us to pass to the limit in each
one of the linear equations of system (\ref{linear Ito}). A posteriori, from
these integral equations, it follows that there is a modification such that
all components are continuous. The proof of existence is complete.
\end{proof}

\section{Girsanov transform\label{sect Girsanov}}

The idea is that equation (\ref{eq Ito}) written in the form
\begin{align*}
dX_{n}  &  =\sigma k_{n-1}X_{n-1}\left(  \frac{1}{\sigma}X_{n-1}%
dt+dW_{n-1}\right)  -\sigma k_{n}X_{n+1}\left(  \frac{1}{\sigma}X_{n}%
dt+dW_{n}\right) \\
&  -\frac{\sigma^{2}}{2}\left(  k_{n}^{2}+k_{n-1}^{2}\right)  X_{n}dt
\end{align*}
becomes equation (\ref{linear Ito}) because the processes $B_{n}\left(
t\right)  :=\frac{1}{\sigma}\int_{0}^{t}X_{n}\left(  s\right)  ds+W_{n}\left(
t\right)  $ are Brownian motions with respect to a new measure $Q$ on $\left(
\Omega,F_{T}\right)  $; and conversely, so both weak existence and weak
uniqueness statements transfer from equation (\ref{linear Ito}) to equation
(\ref{eq Ito}). Equation (\ref{linear Ito}) was also proved to be
\textit{strongly} well posed, but the same problem for the nonlinear model
(\ref{eq Ito}) is open.

Let us give the details. We use results about Girsanov theorem that can be
found in \cite{RevuzYor}, Chapter VIII, and an infinite dimensional version
proved in \cite{Kozlov}, \cite{DZ}.

\subsection{Proof of Theorem \ref{main uniqueness thm}}

Let us prepare the proof with a few remarks. Assume that $\left(
X_{n}\right)  _{n\geq1}$ is an exponentially integrable solution. Since in
particular $E\left[  \int_{0}^{T}\sum_{n=1}^{\infty}X_{n}^{2}\left(  s\right)
ds\right]  <\infty$, the process $L_{t}:=-\frac{1}{\sigma}\sum_{n=1}^{\infty
}\int_{0}^{t}X_{n}\left(  s\right)  dW_{n}\left(  s\right)  $ is well defined,
is a martingale and its quadratic variation $\left[  L,L\right]  _{t}$ is
$\frac{1}{\sigma^{2}}\int_{0}^{t}\sum_{n=1}^{\infty}X_{n}^{2}\left(  s\right)
ds$. Since $E\left[  e^{\frac{1}{2\sigma^{2}}\int_{0}^{T}\sum_{n=1}^{\infty
}X_{n}^{2}\left(  t\right)  dt}\right]  <\infty$, Novikov criterium applies,
so $\exp\left(  L_{t}-\frac{1}{2}\left[  L,L\right]  _{t}\right)  $ is a
strictly positive martingale. Define the probability measure $Q$ on $F_{T}$ by
setting
\begin{equation}
\frac{dQ}{dP}=\exp\left(  L_{T}-\frac{1}{2}\left[  L,L\right]  _{T}\right)  .
\label{QT}%
\end{equation}
Notice also that $Q$ and $P$ are equivalent on $F_{T}$, by the strict
positivity and
\begin{equation}
\frac{dP}{dQ}=\exp\left(  Z_{T}-\frac{1}{2}\left[  Z,Z\right]  _{T}\right)
\label{PT}%
\end{equation}
where
\[
Z_{t}=\sum_{n=1}^{\infty}\int_{0}^{t}\frac{1}{\sigma}X_{n}\left(  s\right)
dB_{n}\left(  s\right)
\]%
\[
B_{n}\left(  t\right)  =W_{n}\left(  t\right)  +\int_{0}^{t}\frac{1}{\sigma
}X_{n}\left(  s\right)  ds.
\]
Indeed $\frac{dP}{dQ}=\exp\left(  -L_{T}+\frac{1}{2}\left[  L,L\right]
_{T}\right)  $ and one can check that $-L_{T}+\frac{1}{2}\left[  L,L\right]
_{T}=Z_{T}-\frac{1}{2}\left[  Z,Z\right]  _{T}$.

Under $Q$, $\left(  B_{n}\left(  t\right)  \right)  _{n\geq1,t\in\left[
0,T\right]  }$ is a sequence of independent Brownian motions. Since
\begin{align*}
\int_{0}^{t}k_{n-1}X_{n-1}\left(  s\right)  dB_{n-1}\left(  s\right)   &
=\int_{0}^{t}k_{n-1}X_{n-1}\left(  s\right)  dW_{n-1}\left(  s\right) \\
&  +\int_{0}^{t}k_{n-1}X_{n-1}\left(  s\right)  X_{n-1}\left(  s\right)  ds
\end{align*}
and similarly for $\int_{0}^{t}k_{n}X_{n+1}\left(  s\right)  dB_{n}\left(
s\right)  $, we see that
\begin{align*}
X_{n}\left(  t\right)   &  =X_{n}\left(  0\right)  +\int_{0}^{t}k_{n-1}%
X_{n-1}\left(  s\right)  dB_{n-1}\left(  s\right)  -\int_{0}^{t}k_{n}%
X_{n+1}\left(  s\right)  dB_{n}\left(  s\right) \\
&  -\int_{0}^{t}\frac{1}{2}\left(  k_{n}^{2}+k_{n-1}^{2}\right)  X_{n}\left(
s\right)  ds.
\end{align*}
This is equation (\ref{linear Ito}). We have proved the first half of the
following lemma:

\begin{lemma}
\label{teo Girsanov}If $\left(  \Omega,F_{t},P,W,X\right)  $ is an
exponentially integrable solution of the nonlinear equation (\ref{eq
stratonovich}), then it is a solution of the linear equation (\ref{linear
Ito}) where the processes
\[
B_{n}\left(  t\right)  =W_{n}\left(  t\right)  +\int_{0}^{t}\frac{1}{\sigma
}X_{n}\left(  s\right)  ds
\]
are a sequence of independent Brownian motions on $\left(  \Omega
,F_{T},Q\right)  $, $Q$ defined by (\ref{QT}). In addition, the process $X$ on
$\left(  \Omega,F_{T},Q\right)  $ satisfies the assumptions of theorem
\ref{teor uniq linear}.
\end{lemma}

\begin{proof}
It remains to prove that conditions (\ref{forth order}) and (\ref{assumption
limit}) hold true. We have
\begin{align*}
E^{Q}\left[  \int_{0}^{T}X_{n}^{4}\left(  t\right)  dt\right]   &
=E^{P}\left[  \mathcal{E}\left(  L\right)  _{T}\int_{0}^{T}X_{n}^{4}\left(
t\right)  dt\right] \\
&  =E^{P}\left[  \exp\left(  L_{T}-\left[  L,L\right]  _{T}+\frac{1}{2}\left[
L,L\right]  _{T}\right)  \int_{0}^{T}X_{n}^{4}\left(  t\right)  dt\right]
\end{align*}%
\[
\leq E^{P}\left[  \exp\left(  2L_{T}-2\left[  L,L\right]  _{T}\right)
\right]  ^{1/2}E^{P}\left[  \left(  \int_{0}^{T}X_{n}^{4}\left(  t\right)
dt\right)  ^{2}\exp\left[  L,L\right]  _{T}\right]  ^{1/2}.
\]
The second factor is finite by the condition of exponential integrability of
$X$. The term $E^{P}\left[  \exp\left(  2L_{T}-2\left[  L,L\right]
_{T}\right)  \right]  $ is equal to one, by Girsanov theorem applied to the
martingale $2L_{t}$. The proof of condition (\ref{forth order}) is complete.
As to condition (\ref{assumption limit}), it follows from the fact that
$E\left[  \int_{0}^{T}\sum_{n=1}^{\infty}X_{n}^{2}\left(  s\right)  ds\right]
<\infty$, a consequence of exponential integrability of $X$. The proof is complete.
\end{proof}

One may also check that
\[
dX_{n}=\sigma k_{n-1}X_{n-1}\circ dB_{n-1}-\sigma k_{n}X_{n+1}\circ dB_{n}%
\]
so the previous computations could be described at the level of Stratonovich calculus.

Let us now prove weak uniqueness (the proof is now classical). Assume that
$\left(  \Omega^{\left(  i\right)  },F_{t}^{\left(  i\right)  },P^{\left(
i\right)  },W^{\left(  i\right)  },X^{\left(  i\right)  }\right)  $, $i=1,2$,
are two exponentially integrable solutions of equation (\ref{eq stratonovich})
with the same initial condition $x\in l^{2}$. Then
\begin{equation}
dX_{n}^{\left(  i\right)  }=\sigma k_{n-1}X_{n-1}^{\left(  i\right)  }%
dB_{n-1}^{\left(  i\right)  }-\sigma k_{n}X_{n+1}^{\left(  i\right)  }%
dB_{n}^{\left(  i\right)  }-\frac{\sigma^{2}}{2}\left(  k_{n}^{2}+k_{n-1}%
^{2}\right)  X_{n}^{\left(  i\right)  }dt \label{eq lin i-esima}%
\end{equation}
where, for each $i=1,2$,
\[
B_{n}^{\left(  i\right)  }\left(  t\right)  =W_{n}^{\left(  i\right)  }\left(
t\right)  +\int_{0}^{t}\frac{1}{\sigma}X_{n}^{\left(  i\right)  }\left(
s\right)  ds
\]
is a sequence of independent Brownian motions on $\left(  \Omega^{\left(
i\right)  },F_{T}^{\left(  i\right)  },Q^{\left(  i\right)  }\right)  $,
$Q^{\left(  i\right)  }$ defined by (\ref{QT}) with respect to $\left(
P^{\left(  i\right)  },W^{\left(  i\right)  },X^{\left(  i\right)  }\right)  $.

We have proved in Theorem \ref{teor uniq linear} that equation (\ref{linear
Ito}) has a unique strong solution. Thus it has uniqueness in law on $C\left(
\left[  0,T\right]  ;\mathbb{R}\right)  ^{\mathbb{N}}$, by Yamada-Watanabe
theorem (see \cite{RevuzYor}, \cite{PrevRoeckner}), namely the laws of
$X^{\left(  i\right)  }$ under $Q^{\left(  i\right)  }$ are the same. The
proof of Yamada-Watanabe theorem in this infinite dimensional context, with
the laws on $C\left(  \left[  0,T\right]  ;\mathbb{R}\right)  ^{\mathbb{N}}$,
is step by step identical to the finite dimensional proof, for instance of
\cite{RevuzYor}, Chapter 9, lemma 1.6 and theorem 1.7. We do not repeat it here.

Given $n\in\mathbb{N}$, $t_{1},...,t_{n}\in\left[  0,T\right]  $ and a
measurable bounded function $f:\left(  l^{2}\right)  ^{n}\rightarrow\mathbb{R}
$, from (\ref{PT}) we have
\begin{align*}
&  E^{P^{\left(  i\right)  }}\left[  f\left(  X^{\left(  i\right)  }\left(
t_{1}\right)  ,...,X^{\left(  i\right)  }\left(  t_{n}\right)  \right)
\right] \\
&  =E^{Q^{\left(  i\right)  }}\left[  \exp\left(  Z_{t}^{\left(  i\right)
}-\frac{1}{2}\left[  Z^{\left(  i\right)  },Z^{\left(  i\right)  }\right]
_{t}\right)  f\left(  X^{\left(  i\right)  }\left(  t_{1}\right)
,...,X^{\left(  i\right)  }\left(  t_{n}\right)  \right)  \right]
\end{align*}
where $Z_{t}^{\left(  i\right)  }:=\sum_{n=1}^{\infty}\int_{0}^{t}\frac
{1}{\sigma}X_{n}^{\left(  i\right)  }\left(  s\right)  dB_{n}^{\left(
i\right)  }\left(  s\right)  $. Under $Q^{\left(  i\right)  }$, the law of
$\left(  Z^{\left(  i\right)  },X^{\left(  i\right)  }\right)  $ on $C\left(
\left[  0,T\right]  ;\mathbb{R}\right)  ^{\mathbb{N}}\times C\left(  \left[
0,T\right]  ;\mathbb{R}\right)  ^{\mathbb{N}}$ is independent of $i=1,2$. A
way to explain this fact is to consider the enlarged system of stochastic
equations made of equation (\ref{eq lin i-esima}) and equation
\[
dZ^{\left(  i\right)  }=\sum_{n=1}^{\infty}\frac{1}{\sigma}X_{n}^{\left(
i\right)  }dB_{n}^{\left(  i\right)  }.
\]
This enlarged system has strong uniqueness, for trivial reasons, and thus also
weak uniqueness by Yamada-Watanabe theorem.

Hence
\[
E^{P^{\left(  1\right)  }}\left[  f\left(  X^{\left(  1\right)  }\left(
t_{1}\right)  ,...,X^{\left(  1\right)  }\left(  t_{n}\right)  \right)
\right]  =E^{P^{\left(  2\right)  }}\left[  f\left(  X^{\left(  2\right)
}\left(  t_{1}\right)  ,...,X^{\left(  2\right)  }\left(  t_{n}\right)
\right)  \right]  .
\]
Thus we have uniqueness of the laws of $X^{\left(  i\right)  }$ on $C\left(
\left[  0,T\right]  ;\mathbb{R}\right)  ^{\mathbb{N}}$. The proof of
uniqueness is complete.

\subsection{Proof of Theorem \ref{teo existence}}

Let $\left(  \Omega,F_{t},Q,B,X\right)  $ be a solution in $L^{\infty}\left(
\Omega\times\left[  0,T\right]  ;l^{2}\right)  $ of the linear equation
(\ref{linear Ito}), provided by Theorem \ref{teo existence linear}. Let us
argue as in the previous subsection but from $Q$ to $P$, namely by introducing
the new measure $P$ on $\left(  \Omega,F_{T}\right)  $ defined as $\frac
{dP}{dQ}=\exp\left(  Z_{T}-\frac{1}{2}\left[  Z,Z\right]  _{T}\right)  $ where
$Z_{t}:=\sum_{n=1}^{\infty}\int_{0}^{t}\frac{1}{\sigma}X_{n}\left(  s\right)
dB_{n}\left(  s\right)  $. Under $P$, the processes
\[
W_{n}\left(  t\right)  :=B_{n}\left(  t\right)  -\int_{0}^{t}\frac{1}{\sigma
}X_{n}\left(  s\right)  ds
\]
are a sequence of independent Brownian motions. We obtain that $\left(
\Omega,F_{t},P,W,X\right)  $ is an $L^{\infty}$-solution of the nonlinear
equation (\ref{eq Ito}). The $L^{\infty}$-property is preserved since $P$ and
$Q$ are equivalent. The proof of existence is complete.

\bibliographystyle{amsplain}

\end{document}